\documentclass{iucrjournals}

\usepackage[utf8]{inputenc}
\usepackage{graphicx}
\usepackage{authblk}
\usepackage{amsmath}
\usepackage{amsfonts}
\usepackage{color}
\usepackage{bm}
\usepackage{tabularx}
\usepackage{url}

\title{Chasing rotational symmetry in the spectre tiling}

\author[1]{Marianne Imperor-Clerc \IUCrCemaillink{marianne.imperor@cnrs.fr}\IUCrOrcidlink{0000-0002-0269-7556}}%
\author[1]{Jean-François Sadoc\IUCrOrcidlink{0000-0002-3332-244X}}%
\author[ ]{Simon Tatham\IUCrCemaillink{anakin@pobox.com}\IUCrOrcidlink{0009-0006-8817-9018}}%

\affil[1]{Laboratoire de Physique des Solides, CNRS and Université Paris-Saclay, 91400 Orsay, France}

\begin{document}
\maketitle

\begin{center}
\today    
\end{center}

\par

\begin{abstract}
We identify the rotational symmetry centers in the spectre tiling. We find that both 2-fold and 3-fold symmetry centers are coexisting, with two different types of 3-fold centers. Local and global symmetry operations are distinguished, and characterized in terms of substitution systems. We introduce a decomposition of the spectre tiling in terms of three metatiles named diabolo, 2-arms and 3-arms, valid at any substitution level. We propose geometrical constructions for generating patches of the spectre tiling with 2-fold and 3-fold global rotational symmetry.
\end{abstract}


\section{Introduction}
Since its discovery, the spectre tiling \cite{smith2023chiral} is attracting a lot of attention because it solved the longstanding `Ein Stein' problem. Indeed, the spectre tile is so far the only known monotile with the amazing property to tile the infinite plane only in a quasiperiodic fashion, by repeating itself in rotated copies turned from each other by multiples of the $2\pi/12$ angle. Contrary to the hat tiling \cite{smith2023aperiodic}, the spectre tiling exhibits no mirror operation and is combining together only rotational symmetries \cite{Kaplan_diffraction_2024}. So far, studies about the spectre tiling have focused on combinatorial properties and substitution schemes in order to prove its quasiperiodicity \cite{smith2023chiral,cheritat_2024}. Finite-state transducer algorithms for handling substitution by treating a tile’s hierarchy
of supertiles in a purely combinatorial fashion have been recently introduced \cite{tatham_2025}. But the link of these properties with rotational symmetries has not been explored yet. 
\par
By definition, when applying a rotation operation around a symmetry center, a spectre tile superimposes on another one. Rotational symmetry centers can be either local or global. Around a global symmetry center, the rotation operation is applying to all tiles of the spectre tiling. In the case of a local symmetry center, the rotation applies only to a finite number of tiles around the center and not to other tiles located further apart of the center. When local symmetry centers are easily observed in any patch of the spectre tiling because they are numerous, global symmetry centers need to be visualized on purpose because they are unique in the infinite tiling. Their unicity comes from the fact that if a tiling contains two global symmetry centers at a finite distance from each other, this tiling is periodic. Conversely, in an infinite quasiperodic tiling like the spectre tiling, two different global symmetry centers are necessarily at infinite distance from each other. 
\par
Our purpose is to fully investigate the rotational symmetries of the spectre tiling, both local and global, and how they are related to combinatorial properties using substitution systems. First we identify the different types of local symmetry centers. Then, by observing extended patches of the spectre tiling generated using substitution schemes adapted from the literature \cite{smith2023aperiodic,tatham_2025}, we identify a decomposition into three metatiles, named diabolo, 2-arms and 3-arms. We demonstrate that this decomposition is present at any substitution level. Finally, we propose geometrical constructions for generating patches of the spectre tiling with 2-fold or 3-fold global rotational symmetry.

\section{Local symmetry centers}
\label{sec:local}
\begin{figure*}
\includegraphics[width=0.9\textwidth]{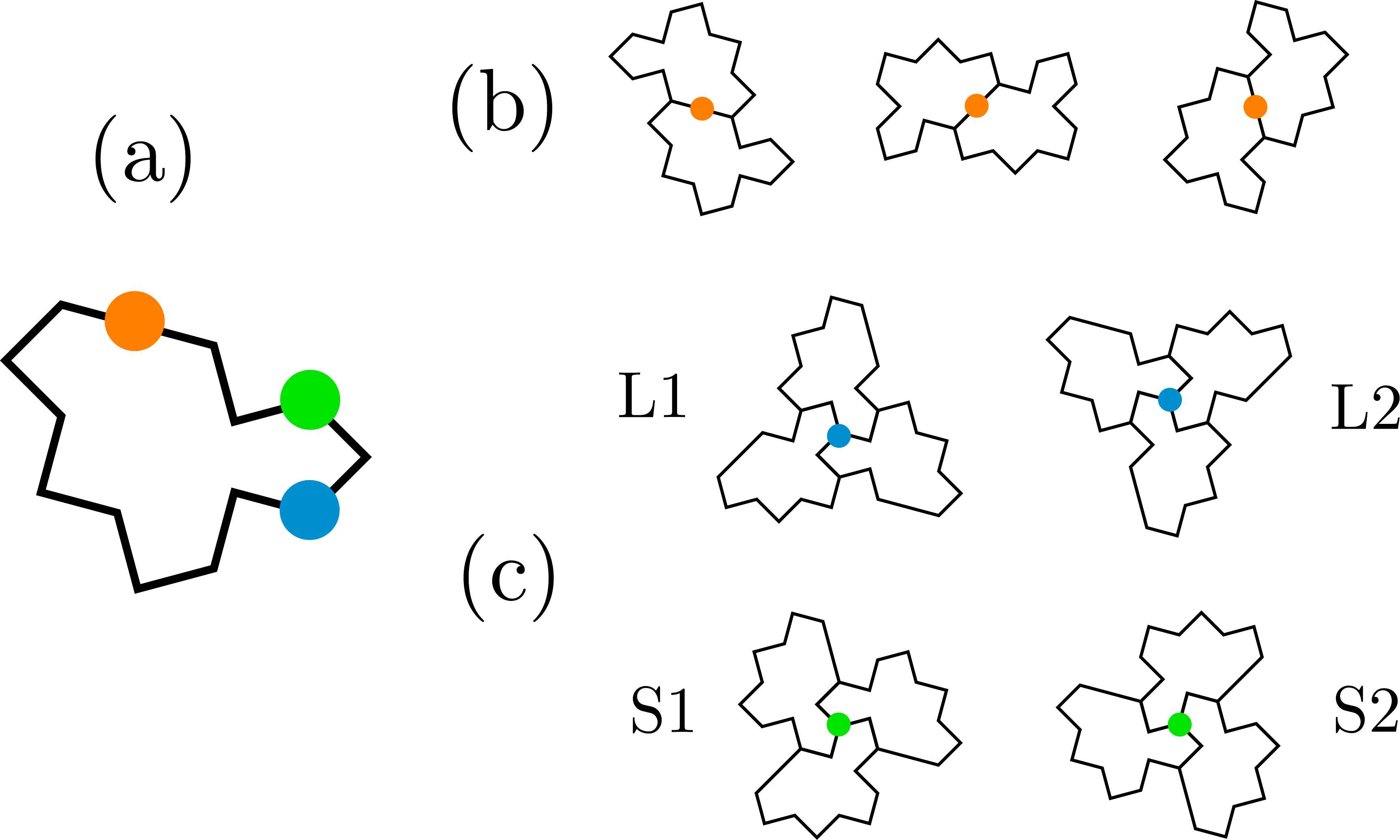}
\centering
\caption{Position of local symmetry centers on the contour of the spectre tile (left chirality, even orientation). (a) The 2-fold center (orange) is located at the middle of the `long' edge of the spectre tile when the two different 3-fold centers (blue and green) have two specific locations. (b) Applying 2-fold rotational symmetry operation around the orange center is giving a set of 2 spectre tiles rotated from each other by $\pi$. In the infinite tiling, this set is taking three possible orientations. (c) Applying 3-fold rotational symmetry operation around a blue or green centers leads to a set of three tiles turned with respect to each other by $2\pi/3$. Note the difference between the two configurations L (blue) and S (green). In the infinite tiling, these sets of three tiles have respectively two possible orientations, (L1, L2) around a blue center and (S1, S2) around a green one.}
\label{fig:monotile}
\end{figure*}

The shape of the spectre monotile is recalled in Figure \ref{fig:monotile}(a). Note that the spectre tile is another name for the Tile(1,1) when referring to other references from the literature \cite{smith2023aperiodic, smith2023chiral, imperor_sadoc_2025}. One recognizes its contour made of a sequence of 14 equal length segments based on the dodecagonal symmetry. One segment is repeated twice, giving rise to one `long' edge in the contour. The spectre monotile is chiral, with two possible handedness. By convention, we assign the left chirality to the shape as shown in Figure \ref{fig:monotile}(a) and all Figures in this paper are depicted with the same left chirality. We identify on the contour of the spectre tile the positions of three different local symmetry centers. The first one is the middle point of the `long' edge (orange point) and it corresponds to a 2-fold center. When 2-fold rotational symmetry is applied to a single spectre tile, a second spectre tile is obtained by rotation of $\pi$ angle around the orange point and the set of two tiles shown in Figure \ref{fig:monotile}(b) is obtained. By construction, this set of two tiles is globally symmetric around the orange point with a rotation angle of $\pi$ (2-fold center). Second, two different points on the contour of the spectre tile (blue and green points) are identified as 3-fold symmetry centers. When a rotation of $\frac{2\pi}{3}$ angle (3-fold symmetry) is applied to a single tile around either a blue or green point, a set of three tiles is obtained as shown in \ref{fig:monotile}(c). Even if these two sets look similar at first sight, a close inspection reveals that they are definitively two different configurations of tiles. The positions on the contour of the spectre tile of the blue and green points are indeed not equivalent, which is somehow related to the chiral nature of the spectre tile itself. To distinguish between these two different local 3-fold centers, we will refer to them as L and S centers for the blue and green points respectively. In the infinite spectre tiling, the L and S local 3-fold centers are present in two possible
orientations, (L1, L2) and (S1, S2), which are rotated from each other by an angle $\pi$ as can be seen in Figure \ref{fig:monotile}(c). 
\par
In Figure \ref{fig:spectre_tiling}(a), some typical region of the spectre tiling is shown. The spectre tile is taking 12 different orientations, obtained by rotation angles of $n 2\pi/12$ where $n$ is an integer. Two subsets of tiles can be distinguished, the ODD and EVEN tiles \cite{smith2023chiral,cheritat_2024,Baake_Mazac_Gahler_2024,imperor_sadoc_2025}. The ODD tiles (M1, M3, M5, M7, M9 and M11) are depicted in pink color, and are related by rotation angles of $n2\pi/6$. They are isolated in the tiling and surrounded by EVEN tiles (M2, M3, M6, M8, M10, M12) with rotation angles of $\pi/6+n2\pi/6$. The number of even tiles is much larger than the number of odd tiles, in a ratio between 7 and 8, which value is the irrational quantity $4+\sqrt{15} \approx 7.873$ in the infinite spectre tiling \cite{smith2023aperiodic,imperor_sadoc_2025}. Note that the absolute orientation of the left-handed spectre tiles is the same in all Figures in this work, in order to keep a consistent definition of the orientations. 
\par
In order to better visualize the rotational symmetries in the spectre tiling, we introduce  a color code which is adapted from the literature \cite{smith2023chiral,tatham_2025}. This color code is based on the initial 9-tile substitution system. As already mentioned, all odd tiles are depicted in pink color. As a simplification, following \cite{Baake_Mazac_Gahler_2024}, we introduce a metatile Ht made of four tiles, one odd tile (pink color) surrounded by three even tiles (white color) as shown in Figure \ref{fig:spectre_tiling}(c). When referring to the literature, the Ht metatile corresponds to grouping together a mystic metatile with a $\Sigma$ tile and a $\Delta$ tile \cite{smith2023chiral}. This metatile Ht is indeed present for all odd tiles in the spectre tiling. It is larger than the mystic metatile (four tiles instead of two) and it incorporates a larger part of the local environment around each odd spectre tile. Its substitution rule is derived by combining together the substitution rules for the group of a mystic metatile with a $\Sigma$ tile and a $\Delta$ tile and it gives a metatile made of three Ht metatiles and 14 even tiles (see Figure \ref{fig:spectre_tiling}(c)). Along with the Ht metatile, only six even tiles are left, with the following names and associated colors: $\Psi$ (light blue), $\Pi$ (dark blue), $\Phi$ (green), $\Theta$ (red), $\Xi$ (orange), and $\Lambda$ (purple). The substitution rule transforms each of these six even tiles into an $\Omega$ metatile (\cite{imperor_sadoc_2025}), but with different coloring for the even tiles (see Figure \ref{fig:spectre_tiling}(c)). Finally, the same left chirality of the spectre tiles is kept at all steps of the substitution scheme (see Figure \ref{fig:spectre_tiling}(d)). This is obtained by generating first the substitution using its original version \cite{smith2023chiral} with alternate handedness at each step. Then the same mirror operation is applied at substitution steps 0, 2, 4, etc ... to get rid of the alternate handedness. In this way, the substitution scheme is completely homochiral, allowing to follow the position of a symmetry center directly for all steps of the substitution, a crucial point in our approach. 
\par
We can now observe where all the local symmetry centers are located in the spectre tiling, by plotting together on the tiling all the orange, blue and green points identified previously. First of all, as can be seen in Figure \ref{fig:spectre_tiling}(a), all local symmetry centers are never located on the contour of odd tiles, Ht metatiles and $\Lambda$ even tiles. Instead, all local symmetry centers are applying between adjacent even tiles. Moreover, there is a strong correlation with the color code from the substitution scheme and the local symmetry centers (see Figure \ref{fig:spectre_tiling}(c)). Sets of three tiles around each L1/L2 3-fold center are made of $\Psi$ (light blue) and $\Pi$ (dark blue) tiles, with the difference that a $\Psi$ tile is carrying a 2-fold center when a $\Pi$ tile do not. All S1/S2 3-fold centers are surrounded by sets of three $\Phi$ (green) tiles. In addition, a subset of isolated $\Phi$ tiles carrying only a local 2-fold center (orange point) is present in the spectre tiling. As a result, the two blue colors (all $\Psi$ and $\Pi$ tiles) and the green color (subset of $\Phi$ tiles) are respectively matching the L and S 3-fold symmetry centers with the sets of three tiles introduced in Figure \ref{fig:monotile}(c). A particularly striking and symmetric configuration is a L set of three $\Psi$ (light blue) tiles surrounded by three S sets of $\Phi$ (green) tiles, as discussed later on in section \ref{sec:global} (see Figure \ref{fig:compare-L-and-S}(a)). Finally, a $\Theta$ (red) and a $\Lambda$ (purple) tile are always glued together in a fixed configuration, as already mentioned in the literature \cite{Baake_Mazac_Gahler_2024}. For the sake of completeness, we may mention that a $\Xi$ (orange) tile is always glued to a $\Pi$ (dark blue) tile with the same fixed configuration.
\par
In conclusion, by cross-linking the positions of the local symmetry centers with the color code from the substitution scheme, we identify the following partition of the spectre tiling. On one side, two types of isolated objects, the Ht metatiles (pink and white), and the pairs of $\Theta$ (red) and $\Lambda$ (purple) tiles. And on the other side, a continuous ensemble around these isolated objects made of the blue ($\Psi$ and $\Pi$), green ($\Phi$) and orange ($\Xi$) even tiles, all derived from each other by applying local 2-fold and 3-fold (L and S) symmetry centers.

\section{Diabolo, 2-arms and 3-arms}
\label{sec:patterns}
Taking a step further, we now investigate at an intermediate length scale how the local symmetry centers are placed around the odd tiles in the spectre tiling. The idea behind this is to use the fact that all odd tiles are sharing the same local environment made of neighboring even tiles that can be represented by a cluster C (see Figure \ref{fig:spectre_tiling}(b)). To do that, we introduce a reference triangle (orange color) around each odd tile as depicted in Figure \ref{fig:orange_triangle}(a). Its initial purpose is to help visualizing the environment around each odd tile. As shown in Figure \ref{fig:orange_triangle}(a), the reference triangle is joining the middle of three `long' edges of three even tiles located around the Ht metatile (see Figure \ref{fig:spectre_tiling}(c)). This triangle is regular, and has the remarkable property that its center is located at the center of the `trefoil' part of the odd tile, when the later is decomposed into two shapes, a `trefoil' and a `bow-tie'. This decomposition of the spectre tile is borrowed to the wheel tiling \cite{nissen_1990}, a tiling strongly related to the spectre tiling. Finally, further away from the central odd tile, more even tiles are represented but with dashed edges, as they may overlap in the spectre tiling with other tiles from adjacent C clusters.
\par
A striking fact is that all the local 2-fold symmetry centers (orange points) in the spectre tiling are located on the vertices of the reference triangles. For most of the odd tiles, all the three vertices of the reference triangle are carrying local 2-fold centers. This configuration is depicted in Figure \ref{fig:orange_triangle}(b) and it corresponds to a $\Omega$ metatile. In addition, a less frequent configuration between two adjacent clusters is present and is depicted in Figure \ref{fig:orange_triangle}(c). In that case, two clusters in relative orientations C1 and C6 are placed in a way that two dashed tiles are overlapping between them. In this specific configuration of two clusters, two vertices of the reference triangles (purple and blue points) do not correspond to local 2-fold centers. Instead, around these purple and blue points, there are three tiles which are glued in a specific fashion with no local rotational symmetry. This less frequent configuration is equivalent to the occurrence in the spectre tiling of a $\Gamma$ metatile on which the blue point is placed on. Indeed, the overlap of two dashed tiles is equivalent to the association of a $\Gamma1$ metatile with a $\Omega2$ metatile as shown in Figure \ref{fig:orange_triangle}(c). Another handy visual reference is the usage of the purple color, with all the purple vertices on the reference triangles placed on the `long' edge of the $\Lambda$ (purple) even tiles. We can also notice that the number of $\Lambda$ tiles is equal to the number of $\Theta$ tiles (red color) and also equal to the number of $\Gamma$ metatiles. 
\par
If we now observe how the reference triangles are arranged in the spectre tiling, the result is striking: the reference triangles are connected together in such a way that only three finite size patterns can be formed (see Figure \ref{fig:diabolo_2arms_3arms}). Based on their shape, we named them as diabolo, 2-arms and 3-arms. Let's describe these three patterns in details (see Figure \ref{fig:diabolo_2arms_3arms}). This observation was double checked for several extended patches of the spectre tiling like the one in Figure \ref{fig:diabolo_2arms_3arms}. Moreover, based on the decomposition in terms of $\Omega$ and $\Gamma$ metatiles, we can introduce three new metatiles for the diabolo, 2-arms and 3-arms. The first one, named as diabolo (blue color), is made of two $\Omega$ metatiles turned by $\pi$ from each other (as $\Omega$1 and $\Omega$4), corresponding to a pattern of two reference triangles connected together. The diabolo metatile takes three different orientations in the spectre tiling. Most importantly, the diabolo metatile has 2-fold symmetry around the common vertex of the two reference triangles, and this very point is identified as the 2-fold global symmetry center of the spectre tiling, as discussed in the next section \ref{sec:global}. The 2-arms (yellow color) and 3-arms (purple color) metatiles are larger and correspond respectively to 8 and 19 reference triangles. The 2-arms contains one $\Gamma$ tile and 7 $\Omega$ tiles when the 3-arms contains 3 $\Gamma$ tiles and 16 $\Omega$ tiles, as detailed in Figure \ref{fig:diabolo_2arms_3arms}. Both the 2-arms and 3-arms metatiles are taking six different orientations in the spectre tiling. Indeed, importantly enough, these two metatiles do not have rotational symmetry. For example, if a 2-arms metatile is rotated by $\pi$, the obtained metatile is not superposable with the initial one. In the same way, the 3-arms metatile is not invariant by a rotation of $2\pi/3$. This fact can also be visualized from the orientations of the single $\Theta$ tile (red color) in the 2-arms and of the three $\Theta$ tiles in the 3-arms. 
\par
Based on these observations, we deduce the important result that the diabolo/2-arms/3-arms metatiles are corresponding to a unique decomposition of the spectre tiling. This can be verified computationally in the following way. Firstly, we can infer the diabolo/2-arms/3-arms decomposition at a lower substitution level, called step 0. The metatiles at step 0 are obtained from the metatiles at step 1 by applying the substituion rule `backwards', using the decomposition into $\Omega$ and $\Gamma$ metatiles at step 1. The result is shown in Figure \ref{fig:deflation-2-fold}(a). The diabolo metatile at step 0, or diabolo$_0$, is a set of two even $\Psi$ spectre tiles turned from each other by the $\pi$ angle. The 2-arms$_0$ contains one odd tile and 8 even tiles when the 3-arms$_0$ is made of three odd tiles and 19 even tiles. These three metatiles (diabolo$_0$/2-arms$_0$/3-arms$_0$) can be also uniquely colored in terms of spectre tiles, or using the Ht metatile as well. One can also notice that 2-arms$_0$ and $\Omega$ metatiles are sharing the same composition in terms of odd and even spectre tiles. 
\par
Then it is possible to prove numerically that every complete spectre tiling of the plane has a unique decomposition into diabolo$_0$, 2-arms$_0$ and 3-arms$_0$ metatiles. This can be done with a finite amount of computation by considering the hierarchical address of a tile within the usual 9-tile substitution system of spectre tiles \cite{smith2023chiral}. A finite number of low-order symbols of a tile's address are enough to determine the orientations of all the surrounding tiles within a specified neighborhood. Even if a very high-order -- or even infinite -- supertile boundary appears within that neighborhood, the transducer techniques from \cite{tatham_2025} can still determine all the local tiles based on a bounded number of symbols. So by examining finitely many symbols, it is possible to test a spectre tile to see whether it can occupy a particular position within one of the three metatiles (diabolo$_0$/2-arms$_0$/3-arms$_0$). Testing all positions within all three metatiles allows us to find cases where a tile cannot be classified, or can occupy more than one possible position. Since all of these tests can be performed based on a finitely long sequence of initial address symbols, it is possible to enumerate all such sequences and verify that in every case a unique answer is obtained.
\par
This analysis is complicated slightly by the fact that two overlapping subsets of the 3-arms$_0$ correspond to the 2-arms$_0$, and many pairs of back-to-back spectre tiles resembling the diabolo$_0$ appear in both the 2-arms$_0$ and 3-arms$_0$ metatiles. However, it is possible to disambiguate using the metatile types from the 9-tile substitution system. A central spectre tile in the 3-arms$_0$ is always a $\Theta$ tile (red color) in the 9-tile system, and conversely, wherever a $\Theta$ tile appears it is always surrounded by a 3-arms$_0$. Similarly, pairs of back-to-back tiles are a diabolo$_0$ if and only if they are $\Psi$ tiles (light blue color).
\par
Under just those constraints (a diabolo$_0$ consists of two $\Psi$ spectre tiles back to back, and a 2-arms$_0$ contains no $\Theta$ spectre), three layers of metatile symbols suffice to classify every possible tile address into exactly one position in one of the three clusters, confirming that every spectre tiling admits a disjoint covering by the clusters.
\par
As mentioned above, the substitution rule for the spectre tiling can always be applied in reverse, so that for any complete plane tiling $T$, there exists a tiling $T'$ which gives rise to $T$ under the substitution rule. $T'$ can be decomposed into the diabolo$_0$, 2-arms$_0$ and 3-arms$_0$, because any tiling can be. Reapplying the substitution rule, we obtain a decomposition of the original $T$ into the larger diabolo, 2-arms and 3-arms at step 1. Iterating this process, it follows by induction that a full spectre tiling can also be decomposed into the $n$th-order substitutions of the diabolo, 2-arms and 3-arms, for any $n>0$.

\section{Global symmetry centers}
\label{sec:global}
We can now determine where global symmetry centers are located in the spectre tiling by identifying their position with respect to the diabolo/2-arms/3-arms decomposition. Then our strategy is to build globally symmetric patches of the spectre tiling of increasing size around each type of symmetry center (2-fold and 3-fold) using the 9-tile substitution system.
\par
\begin{figure*}
\includegraphics[width=1.1\textwidth]{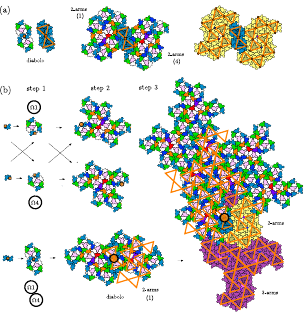}
\centering
\caption{Global 2-fold symmetry center of the spectre tiling. (a) The center of a diabolo is placed at the 2-fold global symmetry point (large orange circle). The diabolo is surrounded by two 2-arms turned from each other by the angle $\pi$ (orientations (1) and (4), see Figure \ref{fig:compare-L-and-S}). (b) Homochiral substitution scheme at steps 1, 2 and 3 for a set of two $\Psi$ tiles. Note how $\Omega$1 and $\Omega$4 orientations alternate at each substitution step. In addition to the diabolo and the two 2-arms patterns, two 3-arms can be identified at step 3.}
  \label{fig:2-fold}
\end{figure*}
\begin{figure*}
\includegraphics[width=1.1\textwidth]{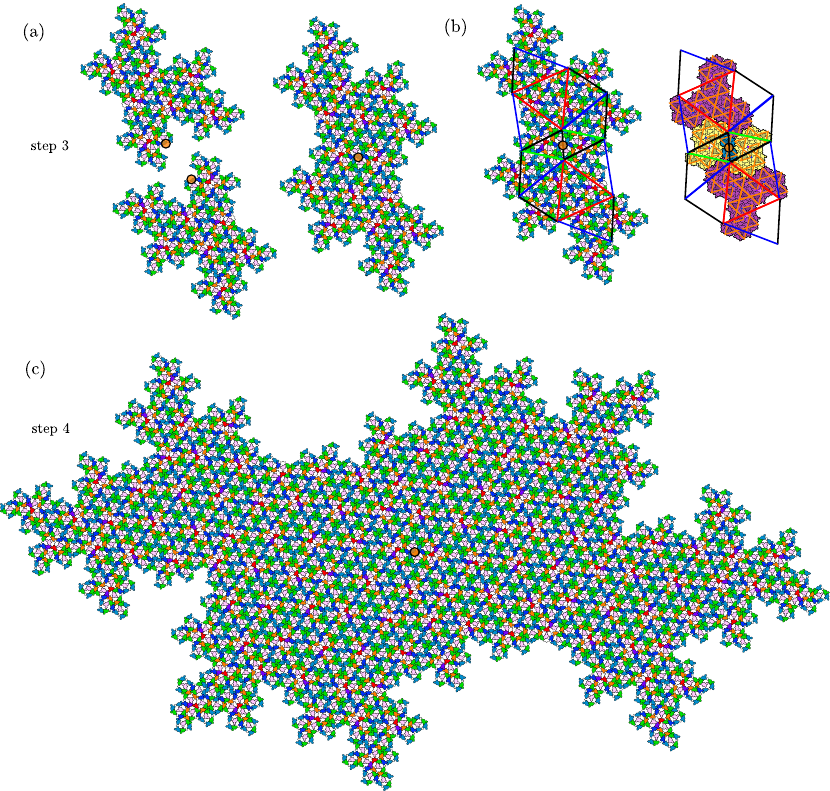}
\centering
\caption{Large patches for global 2-fold symmetry center of the spectre tiling. (a) At step 3, same as in Figure \ref{fig:2-fold}(b), showing how the two patches rotated by $\pi$ angle attach together and form the globally symmetric patch. (b) Same with the superposition of the large triangular network (see Figure \ref{fig:triangular-network}). The rhombus with black and green edges corresponds to the diabolo and the two 2-arms. The large red triangle corresponds to a 3-arm pattern. (c) The globally 2-fold symmetric patch at step 4.}
  \label{fig:2-fold-large}
\end{figure*}

For a global 2-fold symmetry center, the only option for its location is the center of a diabolo metatile (see Figure \ref{fig:2-fold}). As a matter of fact, at the center of a diabolo metatile, there are always two even $\Psi$ tiles turned by the angle $\pi$ around the global 2-fold symmetry center. Starting from the substitution rule for a $\Psi$ tile (see Figure \ref{fig:spectre_tiling}(c)), we can apply it to this set of two $\Psi$ tiles turned from each other by $\pi$ as shown in Figure \ref{fig:2-fold}(b). Note that these two $\Psi$ tiles correspond to a diabolo$_0$ as introduced in section \ref{sec:patterns}. Because the substitution scheme is homochiral, we can identify the location of the 2-fold symmetry center (large orange point) from one substitution step to the next one. The first substitution step shows that the diabolo metatile is surrounded by two 2-arms metatiles. By applying the substitution rule further on, extended patches of the spectre tiling with global 2-fold symmetry can be generated as shown in Figure \ref{fig:2-fold-large}. Finally, the diabolo/2-arms/3-arms decomposition around a global 2-fold symmetry center can be compared between step 0 and step 1 (Figures \ref{fig:deflation-2-fold}(c) and (d)). It is remarkable to observe how the substitution is acting on the diabolo/2-arms/3-arms decomposition producing similar inflated patterns from step 0 to step 1.
\par
Identifying the position of a 3-fold global symmetry center in the spectre tiling is more complicated than for the 2-fold symmetry. As previously mentioned in section \ref{sec:patterns}, the center of a 3-arms metatile is not a 3-fold symmetry center. A possible candidate might be a configuration of three diabolo metatiles rotated from each other by $2\pi/3$, but such configuration of three diabolos is not present in the spectre tiling. Indeed, diabolos can only ever occur between a pair of 2-arms. The only remaining option is then a configuration of three 2-arms metatiles in a 3-fold way as shown in Figure \ref{fig:compare-L-and-S}. Two types of globally 3-fold centers are then identified. First, a L global symmetry center, with a set of three $\Psi$ tiles (light blue) surrounded by three sets of three $\Phi$ tiles (green). In its L1 orientation, it corresponds to three $\Omega$ metatiles ($\Omega$1, $\Omega$3 and $\Omega$5) around the 3-fold symmetry center. The second possibility, a S global symmetry center, is a set of three $\Phi$ tiles (green) surrounded by three sets of three blue tiles ($\Psi$ and $\Pi$). In its S1 orientation, it corresponds to $\Omega$2, $\Omega$4 and $\Omega$6 metatiles. The difference between L and S 3-fold centers is noticeable when observing the orientation of the odd tiles of the $\Omega$ metatiles, which are pointing radially in a L 3-fold center and towards the symmetry center for a S 3-fold center. The way how the three 2-arms are placed differently in the L center (see Figure \ref{fig:compare-L-and-S}(c)) compared to the S center (see Figure \ref{fig:compare-L-and-S}(d)) is also striking, and it confirms our starting hypothesis about the location of 3-fold centers. 
\par
\begin{figure*}
\includegraphics[width=1.1\textwidth]{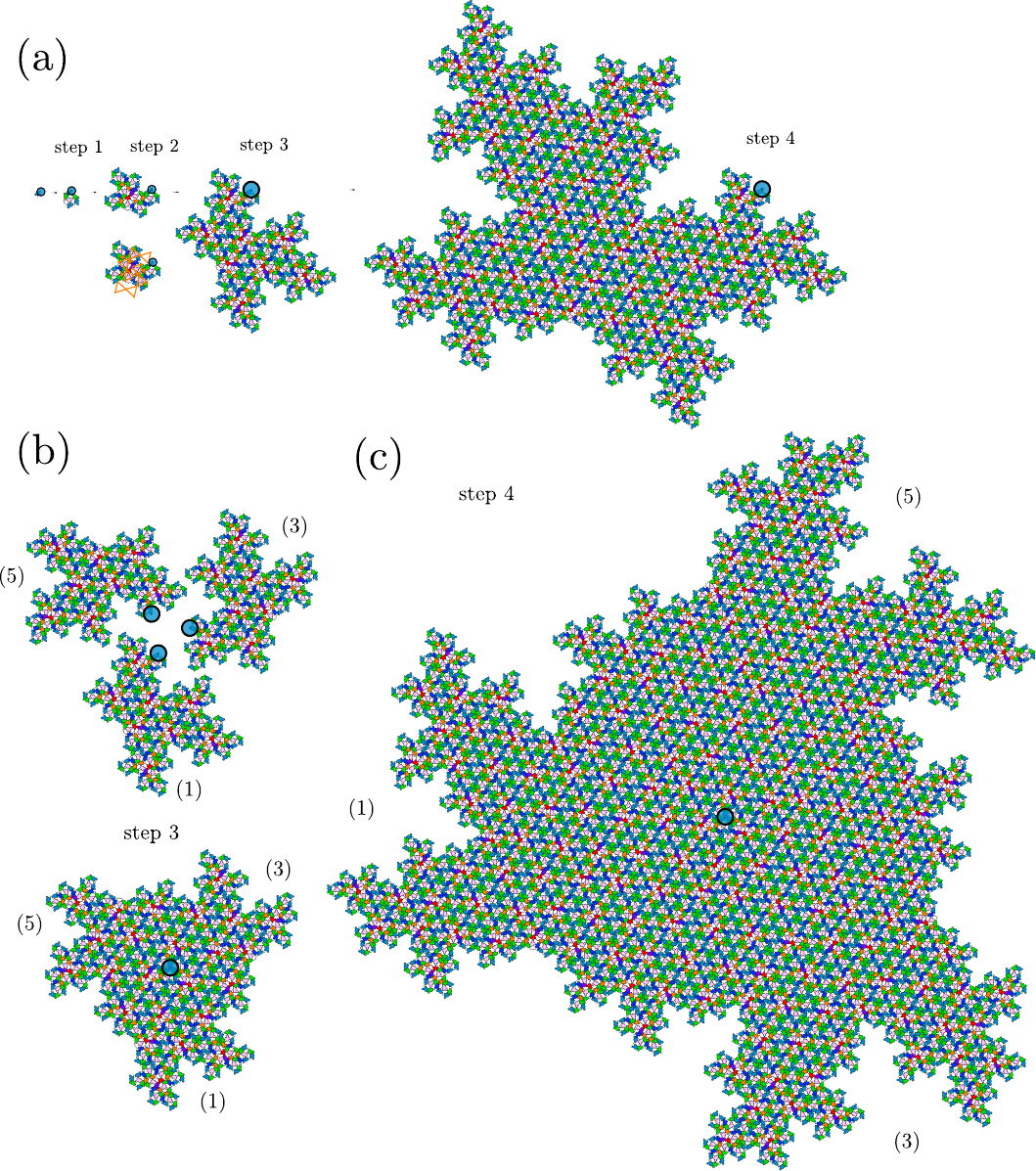}
\centering
\caption{Globally symmetric patches for the L 3-fold type (blue color). (a) Homochiral substitution rule for $\Psi$ tile up to step 4. At each step, the position of the symmetry center (blue point) is uniquely defined as the successive patches are nested from one step to the next one. (b) Globally symmetric 3-fold patch at step 3 obtained by assembling three patches rotated from each other by an angle $2\pi/3$. (c) Similar construction at step 4.}
  \label{fig:3-fold-large}
\end{figure*}

By applying the substitution rules, extended globally symmetric patches with 3-fold symmetry can be build (see Figures \ref{fig:3-fold}, \ref{fig:3-fold-large} and \ref{fig:3-fold-S-green}). For the L 3-fold center, the substitution rule is applied to three $\Psi$ tiles turned from each other by $2\pi/3$ combined together. Note that at each substitution step, the position of the global symmetry center is uniquely determined by identifying a single $\Psi$ tile, in a way that the patch at a previous substitution step is embedded in the patch at the next step. The successive substitution patches fit together like Russian dolls. In this way, the 3-fold L center (large blue circle) is identified on the border of each patch in the substitution scheme of one $\Psi$ tile. For the S 3-fold center, a similar construction is applying using $\Phi$ tiles (green color) instead of $\Psi$ tiles. In that case, the successive patches from one substitution step to the other fit together like Russian dolls but with alternate orientations S1 and S2, as detailed in Figure \ref{fig:3-fold-S-green}. 
\par
Iterating the substitution rule generates a sequence of patches, each globally symmetric. The patch after one substitution is symmetric, and contains the previous patch as a sub-patch at its center of symmetry. Using this as a base case, it follows by induction that the patch after $n$ substitutions is also symmetric, and contains the $(n-1)$th-order patch at its center of symmetry, which in turn contains the $(n-2)$th-order patch in turn, and so on, all the way back to the original patch. This infinite sequence of patches is `consistent', in the sense that if they are all overlaid with their central copy of the original patch aligned, no two patches that cover the same point in the plane disagree about the position of the tile containing that point. Therefore, the union of the whole infinite sequence is a spectre tiling of the plane with the same global symmetry that we began with.
\par
Finally, another geometrical construction can be achieved by joining together all L local 3-fold centers in the spectre tiling. As a result, a triangular network is obtained as illustrated in Figure \ref{fig:triangular-network}. The triangular networks are matching the different global symmetries. This construction helps visualizing the close relationship between local and global symmetry, a notion that could be investigated further on. For example, at step 1, all local 3-fold symmetry centers (L and S) are located at junction between three diabolo/2-arms/3-arms metatiles.  
\par
\begin{figure*}
\includegraphics[width=1.\textwidth]{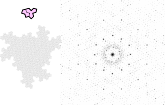}
\centering
\caption{Fourier transform of the spectre tiling using the L1 globally 3-fold symmetric patch (substitution step 3) depicted in Figure \ref{fig:3-fold-large}. A black and white version of the patch including the edges is taken and the Fast Fourier Transform is computed using Image J. Note the orientation of the diffraction pattern with respect to the spectre tiles.}
  \label{fig:FFT}
\end{figure*}

From our analysis, we conclude that the maximum order for global rotational symmetry of the spectre tiling is three. This is in line with the fact that its Fourier transform exhibits 6-fold symmetry as already mentioned in the literature \cite{Kaplan_diffraction_2024,Baake_Gahler_Mazac_diffraction_2025}. Indeed, because Fourier transformation is adding an inversion center, a property known as Friedel's law, 3-fold symmetry in real space gives 6-fold symmetry in Fourier space. As an example, the Fourier transform of the globally L1 3-fold symmetric patch (substitution step 3) is shown in Figure \ref{fig:FFT}. It shows as expected 6-fold symmetry, along with the very signature of chirality, visible from the location of the diffraction peaks.   

\section{Conclusion}
We investigate the rotational symmetry centers of the spectre tiling and find that both 2-fold and 3-fold symmetry centers are coexisting, with two different types (L and S) of 3-fold centers. Our analysis relies on several objects at increasing length scales: local symmetry centers, reference triangle and a decomposition with diabolo/2-arms/3-arms metatiles at different substitution steps. All these objects can be combined together with the substitution scheme to help visualizing the symmetry properties of the spectre tiling. 
\par
\begin{figure*}
\includegraphics[width=1.1\textwidth]{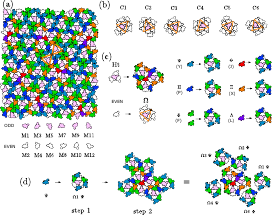}
\centering
\caption{Spectre tiling (left chirality) (a) Finite patch of the spectre tiling. All ODD tiles are depicted in pink and take six possible orientations (M1, M3, M5, M7, M9, M11). All other spectre tiles are in EVEN orientation (M2, M4, M6, M8, M10, M12) and are depicted in different colors depending on the substitution rule. (b) The environment around each odd tile can be represented by a cluster C with six possible orientations (C1, C2, C3, C4, C5, C6). Even tiles with bold edges are always present when even tiles with dashed edges might superimpose with other even tiles from adjacent clusters. The reference triangle is defined in Figure \ref{fig:orange_triangle}(a). (c) Substitution rules for the Ht metatile (one odd tile and three even tiles) into a large metatile and for the six even tiles ($\Psi$, $\Pi$, $\Phi$, $\Theta$, $\Xi$ and $\Lambda$) into a $\Omega$ metatile. Color code for the tiles is detailed in the main text. (d) The two first steps of the substitution for tile $\Psi$ (light blue color). The decomposition of the substitution rule with one Ht metatile and five neighboring tiles is shown on the right.}
\label{fig:spectre_tiling}
\end{figure*}
\begin{figure*}
\includegraphics[width=0.67\textwidth]{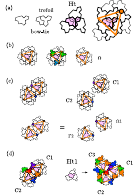}
\centering
\caption{Reference triangle around each odd tile. (a) The reference triangle (orange) is a regular triangle which center is located at the center of the `trefoil' inside the odd tile. Its three vertices are located at the middle of `long' edges of three even tiles. Even tiles with bold edges all belong to the cluster defined by the odd tile in the middle. Other even tiles with dashed edges are present in the spectre tiling but might be attributed to the same cluster or to other neighboring clusters, depending on the configuration. (b) The most common configuration for a cluster is that three local 2-fold centers are located at the three vertices of its reference triangle. These three 2-fold centers are located on a $\Psi$ (light blue) and two $\Phi$ (green) tiles. This configuration also corresponds to a $\Omega$ metatile. (c) The less common configuration is a combination of two adjacent clusters, where two vertices (purple and dark blue) on the reference triangles do not correspond to local 2-fold centers. Note that in between the two odd tiles two even tiles are repeated by translation. In this combination, the two clusters are rotated from each other by $\pi/3$, like C1 and C2, which is equivalent to the association of $\Omega1$ and $\Gamma2$ metatiles. (d) Same configuration as in (c) using the substitution color code (see Figure \ref{fig:spectre_tiling}) and showing the link with the substitution rule for a Ht metatile.} 
\label{fig:orange_triangle}
\end{figure*}
\begin{figure*}
\includegraphics[width=0.9\textwidth]{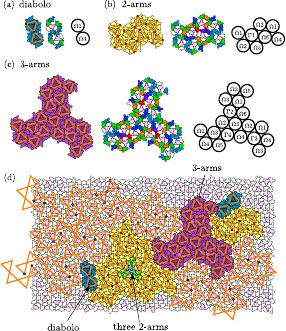}
\centering
\caption{The diabolo/2-arms/3-arms decomposition at substitution step 1: The diabolo metatile (blue) (a), the 2-arms metatile (yellow) (b) and the 3-arms metatile (purple). Patterns formed by the reference triangles are superimposed. The metatiles are also represented with the color code for the spectre tiles (see Figure \ref{fig:spectre_tiling}) and with the decomposition into $\Omega$ and $\Gamma$ metatiles, where the numbers from 1 to 6 are coding their orientations. (d) An extended patch of the spectre tiling showing how all reference triangles connect to each other in the only three possible patterns of the decomposition. All the local 2-fold centers (orange points) are on vertices of reference triangles. The remaining vertices (blue and purple points) do not correspond to 2-fold symmetry centers (see also Figure \ref{fig:orange_triangle}(c)). Blue points are located on the $\Gamma$ metatiles. A configuration of three 2-arms metatiles rotated by $2\pi/3$ around a L1 local 3-fold center is indicated by an arrow (see Figure \ref{fig:compare-L-and-S}(a)).}
  \label{fig:diabolo_2arms_3arms}
\end{figure*}
\begin{figure*}
\includegraphics[width=1.\textwidth]{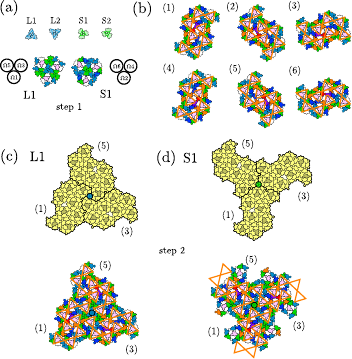}
\centering
\caption{Comparison of the two types (L and S) of 3-fold global symmetry centers. (a) Sets of three spectre tiles for 3-fold symmetry: $\Psi$ tiled for L center (blue) and $\Phi$( tiles for S-center (green). Note the rotation angle of $\pi$/3 between the two possible orientations L1/L2 and S1/S2. Globally symmetric patches at step 1 for L and S centers (in orientations L1 and S1). (b) The six orientations of the 2-arms metatile. Note that it is not 2-fold symmetric, as can be visualized by the orientation of the $\Theta$ tile in red color (see Figure \ref{fig:spectre_tiling}). (c) L 3-fold globally symmetric patch at step 2 is obtained by the combination of three 2-arms turned from each other by $2\pi/3$ (orientations (1), (3) and (5)). (d) S 3-fold globally symmetric patch at step 2 is obtained by combining exactly the same three 2-arms but with a different relative placement.}
  \label{fig:compare-L-and-S}
\end{figure*}
\begin{figure*}
\includegraphics[width=1.05\textwidth]{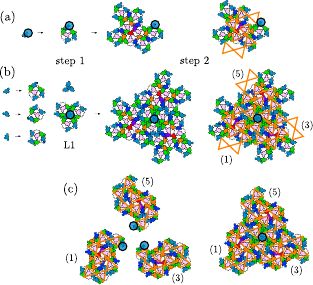}
\centering
\caption{Substitution for the global 3-fold symmetric center of type L. (a) Homochiral substitution rule at steps 1 and 2 for a $\Psi$ tile. (b) Applying the same substitution rule to a set of three $\Psi$ tiles forming a L1 set (see Figure \ref{fig:compare-L-and-S}(a)). (c) Equivalence with the construction using three 2-arms (see Figure \ref{fig:compare-L-and-S}).}
  \label{fig:3-fold}
\end{figure*}
\begin{figure*}
\includegraphics[width=1.1\textwidth]{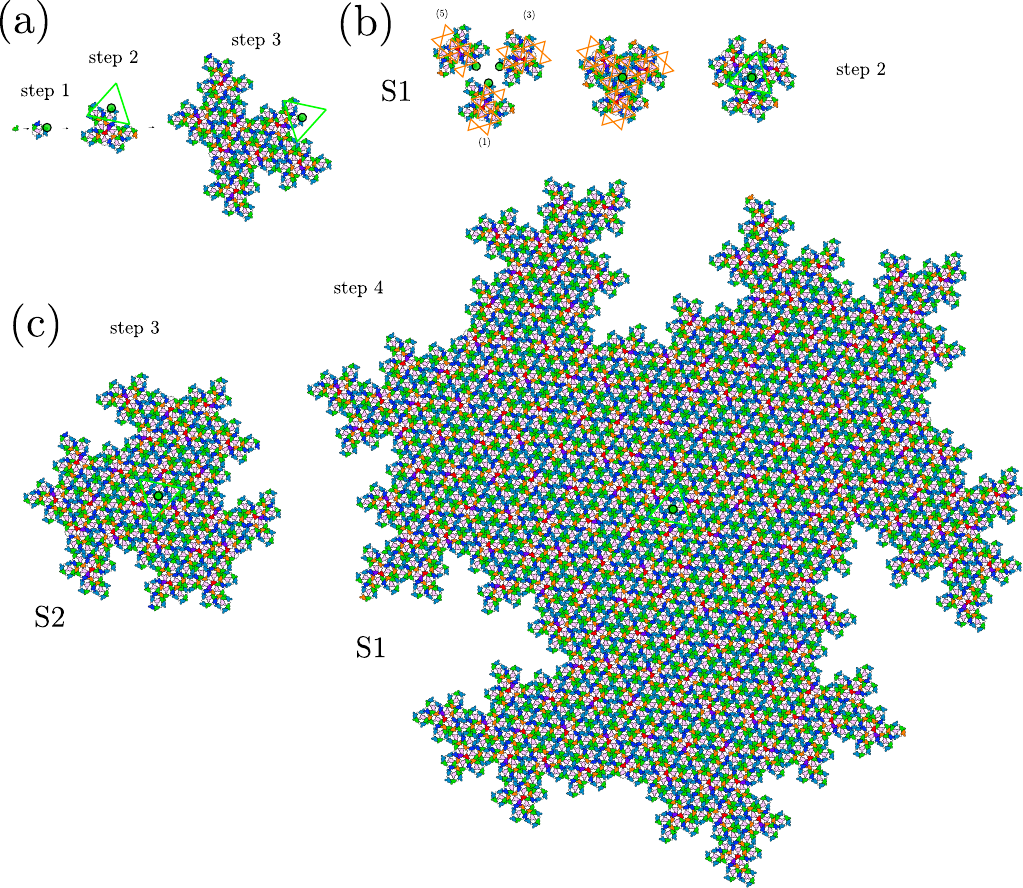}
\centering
\caption{Globally 3-fold symmetry patches for the S type (green color) (a) Homochiral substitution for the $\Phi$ tile. The position of the symmetry center is uniquely defined at each substitution step from the nested patches between two successive steps, but the orientation alternates between S1 and S2 (rotation of $\pi/3$ angle). (b) Construction of the 3-fold symmetric patch at step 2 (see also Figure \ref{fig:compare-L-and-S}). (c) Extended globally symmetric patches at steps 3 and 4. The regular green triangles are centered on the S 3-fold symmetry centers and are part of the triangular network described in Figure \ref{fig:triangular-network}.}
  \label{fig:3-fold-S-green}
\end{figure*}
\begin{figure*}
\includegraphics[width=1.1\textwidth]{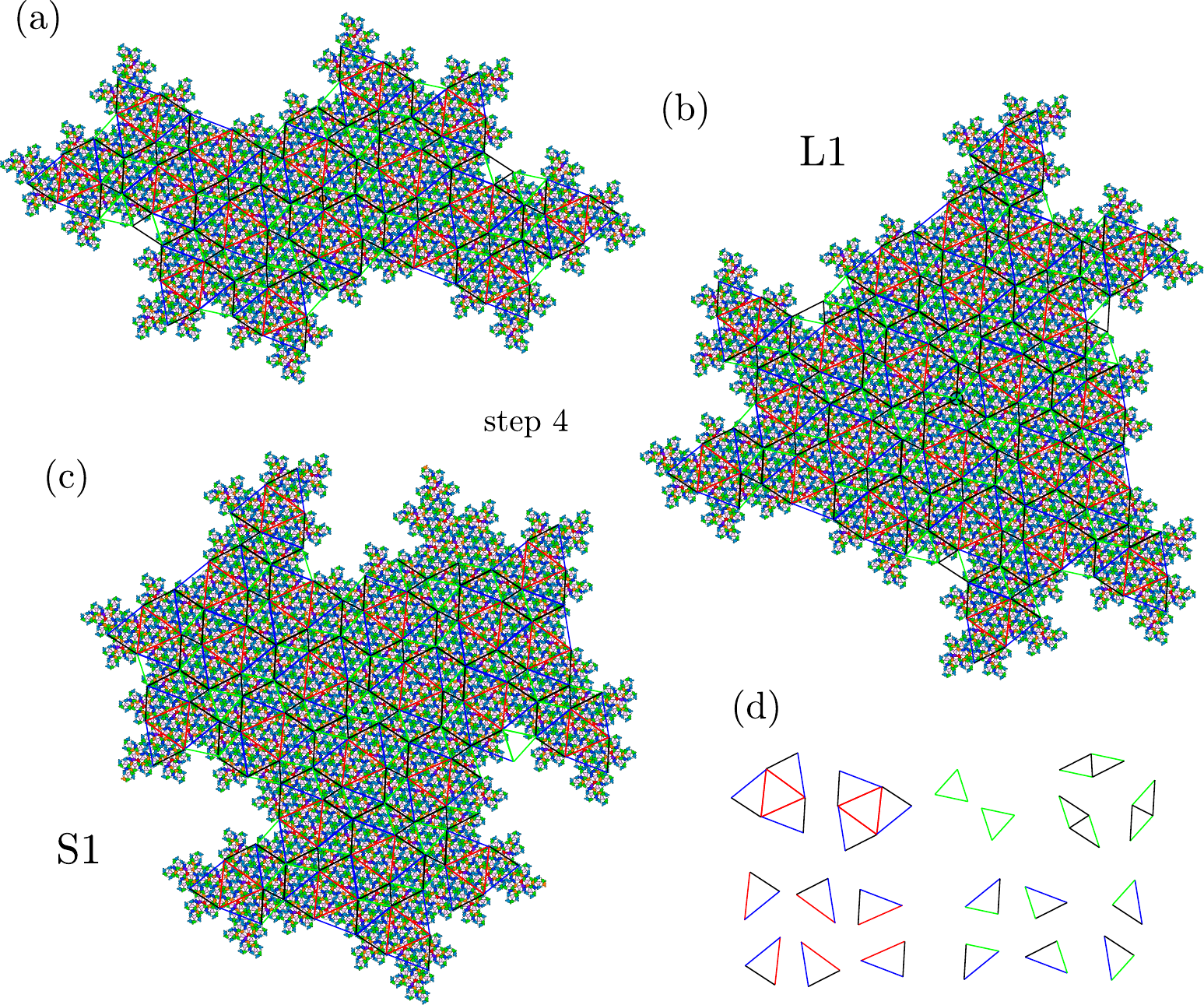}
\centering
\caption{Triangular network joining all L local 3-fold centers. (a) Global 2-fold symmetric patch at step 4 (see Figure \ref{fig:2-fold-large}). (b) Global L 3-fold symmetric patch at step 4. (c) Global S 3-fold symmetric patch at step 4. Different configurations in the triangular networks are observed depending on the symmetry. (d) The finite set of triangles in the network. The five different edge lengths are depicted with colors. The green triangle corresponds to a S-center (see Figure \ref{fig:3-fold-S-green}) when the green/black rhombus corresponds to a diabolo and the red regular triangle to a 3-arms (see Figure \ref{fig:2-fold-large}). } 
  \label{fig:triangular-network}
\end{figure*}

Our approach confirms that different sets of substitution rules can be constructed. For example, the usage of the Ht metatile instead of the Mystic metatile is a new possibility. Moreover, the set of substitution rules presented in Figure \ref{fig:spectre_tiling} might be further optimized in order to highlight the rotational symmetries. For example, it would be better to separate the $\Phi$ (green) tile into two categories, one corresponding to the tiles carrying the 3-fold centers, and another one containing the isolated $\Phi$ tiles, which are all glued in the same way to the $\Pi$ (dark blue) tiles. Another promising perspective is to make usage of the triangular network between local L 3-fold symmetry centers for building other substitution rules but at a much larger length scale than the level of a single spectre tile.  
\par
We propose a general method to construct globally symmetric patches of the spectre tiling around the 2-fold and the L and S 3-fold centers. This construction can be iterated to any substitution step. Maybe other symmetrical configurations for global symmetry are possible, and, in our opinion, the `unicity' question remains open. Homochirality for the substitution rules is an important feature because it allows to identify at all substitution steps identical patches with the same orientation along with the position of the global symmetry center. Finally, our analysis evidences a deep connection between local and global symmetry centers. Here, computation of the average hyperslope of a patch \cite{Imperor_PRB_2024} might be useful to investigate this point further on. Average hyperslope provides four independent real numbers and it is a way to quantify the composition of any patch of the spectre tiling. Comparing the hyperslopes of patches with different global symmetries at successive substitution steps should help determining the symmetry elements in the 4D lift hyperspace.
\par
\begin{figure*}
\includegraphics[width=1.1\textwidth]{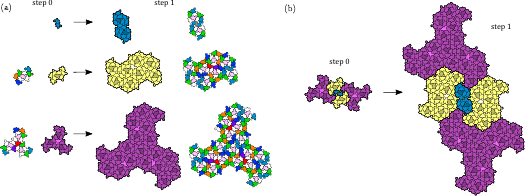}
\includegraphics[width=1.1\textwidth]{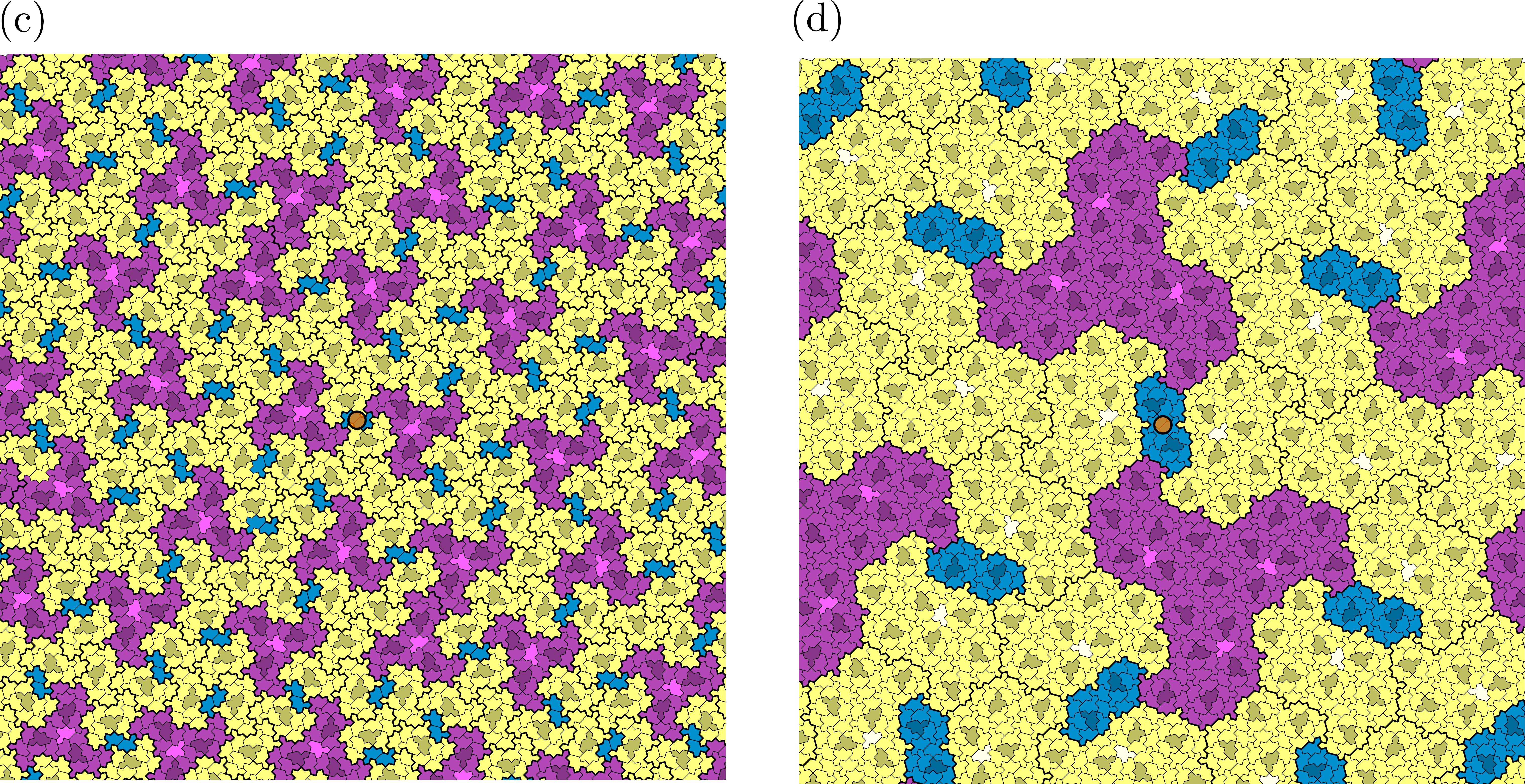}
\centering
\caption{Diabolo/2-arms/3-arms decomposition and 2-fold global symmetry. (a) Decomposition into diabolo/2-arms/3-arms at substitution steps 0 and 1. At step 0, the metatiles are called diabolo$_0$, 2-arms$_0$ and 3-arms$_0$. Note inside the metatiles the distribution of the $\Theta$ tiles, which are initially colored in red, and changed to light yellow in the 2-arms and to light purple in the 3-arms in the blue/yellow/purple color code. (b) Around a global 2-fold symmetry center (orange point), at any substitution step, a diabolo is placed at the symmetry center and is surrounded by two 2-arms and two 3-arms (see Figure \ref{fig:2-fold}). In (c) and (d), the diabolo/2-arms/3-arms decomposition of the same 2-fold globally symmetric patch of the spectre tiling are compared at step 0 (c) and step 1 (d). Note that the patch of spectre tiles is exactly the same in (c) and (d).}
  \label{fig:deflation-2-fold}
\end{figure*}
\begin{figure*}
\includegraphics[width=1.\textwidth]{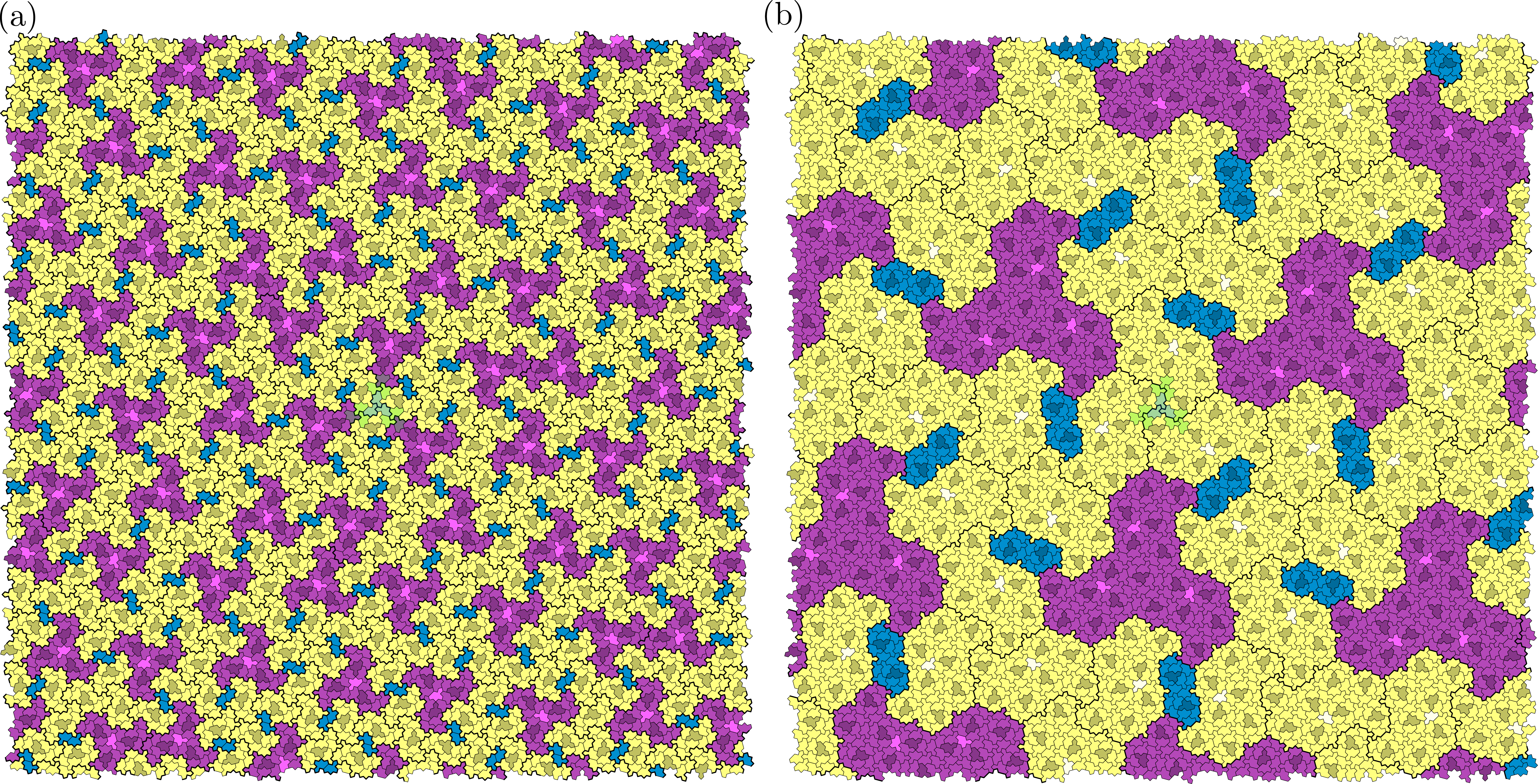}
\includegraphics[width=1.\textwidth]{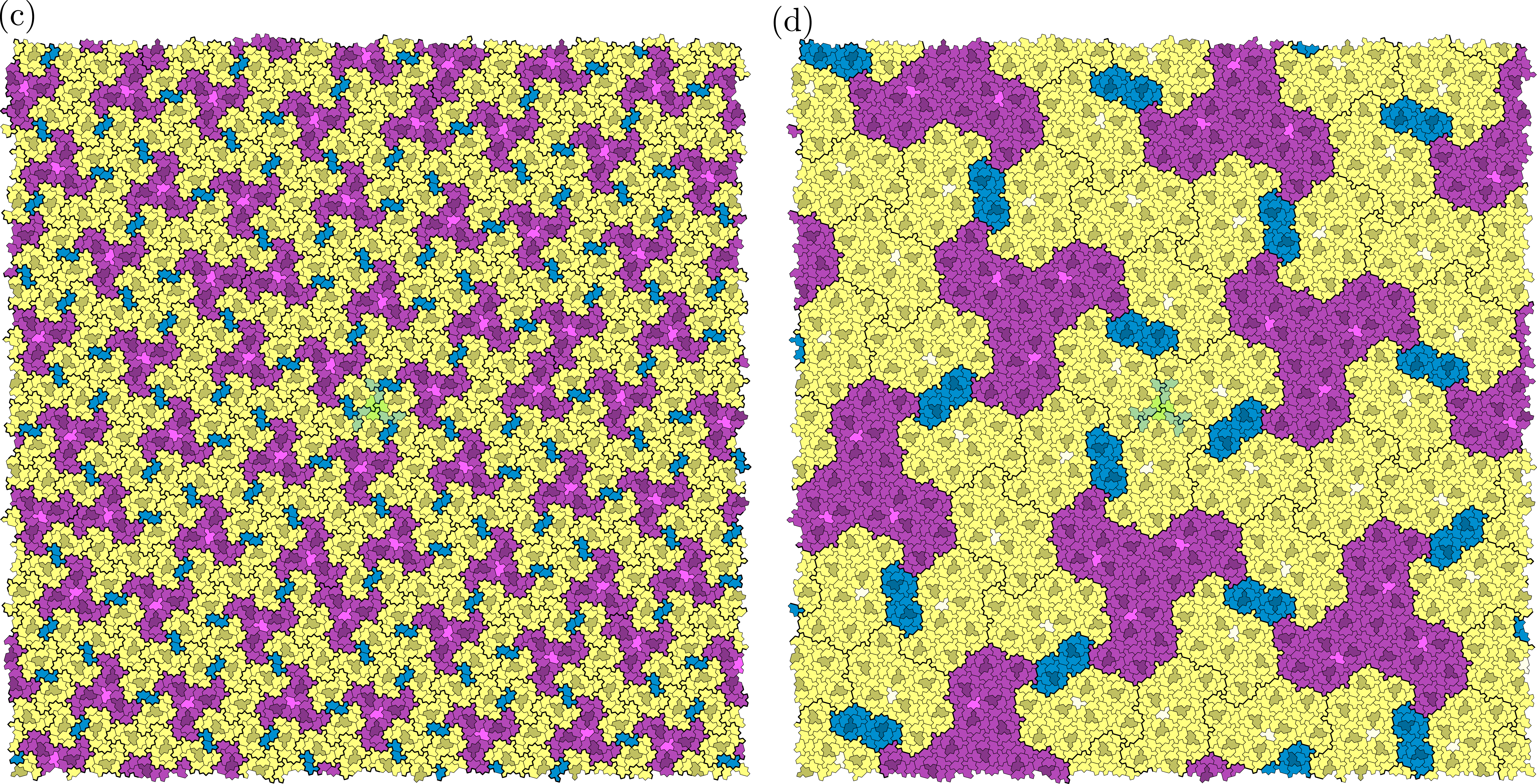}
\centering
\caption{Diabolo/2-arms/3-arms decomposition for 3-fold global symmetry. (a) L 3-fold center at step 0. (b) L 3-fold center at step 1. (a) S 3-fold center at step 0. (b) S 3-fold center at step 1. Note that the spectre tiling is exactly the same in (a)/(b) and (c)/(d). Only the diabolo/2-arms/3-arms decomposition is changing from step 0 to step 1.}
  \label{fig:blue_yellow_purple_3fold}
\end{figure*}

Since any spectre tiling can be decomposed into the diabolo, 2-arms and 3-arms metatiles as shown in Figures \ref{fig:deflation-2-fold} and \ref{fig:blue_yellow_purple_3fold}, it follows by induction that any spectre tiling can also be decomposed into the $n$th-order substitution of those metatiles. Moreover, these three metatiles admit a substitution system in their own right, in which each large metatile deflates to a set of whole small metatiles. A publication demonstrating this is in preparation \cite{Tatham_2026}.
\par
To conclude, the geometrical properties of the spectre tiling, combining chirality and rotational symmetry without periodicity, makes it a very promising candidate as a new class of structure in material science. For example, chiral optical diffraction is a promising application \cite{Moritake_2026}. In this context, we hope that documenting the rotational symmetries of the spectre tiling will help the investigation of this new quasicrystalline materials. 

\begin{acknowledgements}
We are very grateful to Boris K. Horvat  for many stimulating discussions and to Arnaud Chéritat and Anuradha Jagannathan for helpful suggestions for improving this manuscript. Figures are created using the open source software Inkscape. Fast Fourier Transform (FFT) is performed with the open source software ImageJ.
\end{acknowledgements}

\DataAvailability{Python code for generating patches by substitution with the color code is available at \url{https://github.com/marimperorclerc/homochiral}. It is adapted from the open source code available at \url{https://github.com/shrx/spectre}, which algorithm is the original version with alternate chirality between successive substitution steps \cite{smith2023chiral}. The homochiral substitution is derived by applying a vertical mirror at substitution steps 0, 2, 4, etc ... in order to keep the left handedness at all steps. Note that an overall anticlockwise rotation of $\pi/4$ is applied afterwards to all patches to be consistent with the tiles orientation as defined in Figure \ref{fig:spectre_tiling}.
\par
Rust code to verify the correspondence between the diabolo/2-arms/3-arms metatiles decomposition and the existing Spectre tiling substitution system is available at \url{https://github.com/sgtatham/spectre-cluster-classifier}. That program reuses a more general library of finite-state transducer code \cite{tatham_2025}, available from the standard Rust package repository at \href{https://crates.io/}{\texttt{crates.io}} under the name \href{https://crates.io/crates/substitution-tiling-transducers}{\texttt{substitution-tiling-transducers}}.
}

\bibliography{biblio} 


\newpage

\end{document}